\documentclass[12pt]{article}
\usepackage{amsmath}
\usepackage{amssymb}
\usepackage{amsmath,amssymb,amsbsy,amsfonts,amsthm,latexsym,
amsopn,amstext,amsxtra,euscript,amscd}

\newtheorem{lem}{Lemma}
\newtheorem{lemma}[lem]{Lemma}

\newtheorem{thm}{Theorem}
\newtheorem{theorem}[thm]{Theorem}
\newtheorem{cor}{Corollary}

\begin{document}

\title{On small distances between ordinates of zeros of $\zeta(s)$ and $\zeta'(s)$}
\author{{M.~Z.~Garaev}\\
\normalsize{Instituto de Matem{\'a}ticas,  UNAM}
\\
\normalsize{Campus Morelia, Apartado Postal 61-3 (Xangari)}
\\
\normalsize{C.P. 58089, Morelia, Michoac{\'a}n}\\
\normalsize{M{\'e}xico} \\
\normalsize{\tt garaev@matmor.unam.mx} \\
\and\\
{C. Y. Y{\i}ld{\i}r{\i}m}
\\
\normalsize{Department of Mathematics}
\\
\normalsize{Bo\~{g}azi\c{c}i University, Bebek, Istanbul 34342}\\
\normalsize{\& Feza G\"ursey Enstit\"us\"u, \c{C}engelk\"{o}y}\\
\normalsize{ Istanbul, P.K. 6, 81220}\\
\normalsize{Turkey}
\\
\normalsize{\tt yalciny@boun.edu.tr}
 }

\date{\empty}

\pagenumbering{arabic}

\maketitle

\begin{abstract}
We prove that for any zero $\beta'+i\gamma'$ of $\zeta'(s)$ there
exists a zero $\beta+i\gamma$ of $\zeta(s)$ such that
$|\gamma-\gamma'|\ll \sqrt{|\beta'-\tfrac{1}{2}|},$ and we provide
some other related results.
\end{abstract}

\paragraph*{2000 Mathematics Subject Classification:}  11M26 (11M06).

\paragraph*{Key words:} Riemann zeta-function, zeros of $\zeta(s)$ and $\zeta'(s).$

\section{Introduction}

In this paper $s=\sigma+it$ will denote a complex variable, where
$\sigma$ and $t$ are real, and $T$ will denote a large parameter.

The relations between the zeros of a function and the zeros of its
derivatives have been the object of much study. The case of the
Riemann zeta-function $\zeta(s)$ presents many puzzles beginning
with the Riemann hypothesis (RH). Speiser~\cite{Spei} showed that RH
is equivalent to $\zeta'(s)$ having no zeros in $0<\sigma<
\frac{1}{2}.$  From Riemann's original work (proofs for some parts
of which were provided later by other mathematicians), it is
well-known that the non-trivial zeros of $\zeta(s),$ to be denoted
by $\rho=\beta+i\gamma,$ are to be found only in the critical strip,
i.e. $0\le\beta\le 1,$ and the number of non-trivial zeros with
$\gamma\in [0,T]$ is
$$
N(T)=\frac{T}{2\pi}\log\frac{T}{2\pi e}+ \frac{7}{8}+ S(T)
+O(\frac{1}{T})
$$
as $T\to\infty$. Here for $t$ not the ordinate of a zero, $S(t)
:=\frac{1}{\pi}\arg\zeta(\frac{1}{2}+it)$ obtained by continuous
variation along the line segments joining $2,\, 2+it, \,
\frac{1}{2}+it,$ starting with the value $0$; if $t$ is the ordinate
of a zeta zero, $S(t) :=S(t+0)$. It is also well-known that
$S(T)=O(\log T)$. Titchmarsh~\cite[Theorem 11.5 (C)]{T} established
the existence of a constant $E$, between $2$ and $3$, such that
$\zeta'(s)$ does not vanish in the half-plane $\sigma>E,$ while
$\zeta'(s)$ has infinitely many zeros in any strip between
$\sigma=1$ and $\sigma=E.$ Berndt~\cite{Bernd} showed that the
number of non-real zeros of $\zeta'(s),$ which are to be denoted by
$\rho'=\beta'+i\gamma',$ with $\gamma\in[0,T]$ is
$$
N'(T)=\frac{T}{2\pi}\log\frac{T}{4\pi e}+O(\log T).
$$
Levinson and Montgomery~\cite{LM} in addition to proving a
quantified version of Speiser's theorem and that the only zeros of
$\zeta'(s)$ in $\sigma\le 0$ are its `trivial zeros' on the negative
real axis which occur between the trivial zeros of $\zeta(s),$
obtained results revealing that the zeros of $\zeta'(s)$ are mostly
clustered around $\sigma=\frac{1}{2},$ and most of the non-real
zeros of $\zeta'(s)$ lie to the right of
$\sigma=\frac{1}{2}-\frac{w(t)}{\log t},$ where $w(t)\to\infty$ as
$t\to\infty.$ From the fact that $\Re\frac{\zeta'(s)}{\zeta(s)}<0$
on $\sigma=\frac{1}{2},$ except at zeros of $\zeta(s),$ they
observed that $\zeta'(\frac{1}{2}+i\gamma')=0$ can occur only if
$\frac{1}{2}+i\gamma'$ is a multiple zero of $\zeta(s).$ Levinson
and Montgomery also proved
$$
\sum_{0<\gamma'\le T}(\beta'-\tfrac{1}{2})\sim
\frac{T}{2\pi}\log\log T,
$$
which has the immediate interpretation that $\beta'-\frac{1}{2}$ is
often much larger than the average gap between the consecutive zeros
of $\zeta(s).$ In~\cite{CG} Conrey and Ghosh showed that for any
fixed $\nu>0,$ a positive proportion of zeros of $\zeta'(s)$ are in
the region $\sigma\ge \frac{1}{2}+\frac{\nu}{\log t}.$ We note that
the works cited above (except for Titchmarsh's book) deal more
generally with $\zeta^{(k)}(s)$ and contain other results which we
have not mentioned here.

Soundararajan~\cite{Sound} addressed these matters expressing his
belief that the magnitude of $\beta'-\frac{1}{2}$ is usually of
order $\frac{1}{\log \gamma'},$ and the average is high because of
few zeros which are abnormally distant from $\sigma=\frac{1}{2}.$ He
also wrote to the effect that, the more distant $\rho'$ is from the
critical line the larger the gap between the two zeros of $\zeta(s)$
which straddle $\rho'.$ Soundararajan announced two conjectures:

\bigskip

{\bf Conjecture A.} {\it For $\nu\in \mathbb{R},$ let
$$
m^{-}(\nu)=\liminf_{T\to\infty}\frac{1}{N'(T)}\#\{\rho': \,
\beta'\le \frac{1}{2}+\frac{\nu}{\log T}, \quad 0\le \gamma'\le T\}
$$
and $m^{+}(\nu)$ by replacing $\liminf$ by $\limsup$ in the above.
Then  for all $\nu$ we have $m^{-}(\nu)=m^{+}(\nu)=:m(\nu).$
Further, $m(\nu)$ is a nonnegative, nondecreasing, continuous
function with the properties: $m(\nu)=0$ for $\nu\le 0$,
$0<m(\nu)<1$ for $\nu>0,$ and $m(\nu)\to 1$ as $\nu\to \infty.$}

\bigskip

{\bf Conjecture B.} {\it Assume RH. The following two statements are
equivalent:
\begin{itemize}
\item[\rm (i)] \quad $\liminf\limits_{\gamma'\to\infty}
(\beta'-\tfrac{1}{2})\log\gamma'=0;$
\item[\rm (ii)] \quad $\liminf\limits_{\gamma\to\infty}
(\gamma^+-\gamma)\log\gamma=0,$ where $\gamma^+$ is the least
ordinate of a zero of $\zeta(s)$ with $\gamma^+>\gamma.$
\end{itemize}}

Towards these conjectures he showed that there exists a constant $C$
such that $m^{-}(C)>0$ unless RH is `badly violated', and assuming
RH he obtained $m^{-}(\nu)>0$ for $\nu\ge 2.6.$ Zhang~\cite{Zh} made
considerable progress for Conjecture A by proving unconditionally
that $m^{-}(\nu)>0$ for sufficiently large $\nu.$ Assuming RH and
Montgomery's ~\cite{Monty} pair correlation conjecture in the weak
form
$$
\liminf\limits_{T\to\infty}
\frac{1}{N(T)}\sum_{\substack{\gamma_n\le T\\
\gamma_{n+1}-\gamma_n\le \frac{\alpha}{\log T}}}1 >0
$$
for any fixed $\alpha>0,$ Zhang also showed that $m^{-}(\nu)>0$ for
any $\nu>0,$ and Feng~\cite{Feng} was able to dispense with the
assumption of RH in obtaining this result. Here and in what follows
we use the notation that the non-trivial zeros
$\rho_n=\beta_n+i\gamma_n$ of $\zeta(s)$ in the upper half-plane are
indexed as $ 0<\gamma_1\le \gamma_2\le\ldots,$ with the
understanding that the ordinate of a zero of multiplicity $m$
appears $m$ times consecutively in this sequence. Moreover,
Zhang~\cite{Zh} showed under RH that when $\alpha_1$ and $\alpha_2$
are positive constants satisfying $\alpha_1<2\pi$ and
$\alpha_2>\alpha_1\left(1-\sqrt{\frac{\alpha_1}{2\pi}}\right)^{-1},$
if it happens that $(\gamma^+-\gamma)\log\gamma<\alpha_1 $ for
$\rho$ with sufficiently large $\gamma,$ then there exists $\rho'$
such that $|\rho'-\rho|<\alpha_2(\log\gamma)^{-1},$ thereby proving
that ``{\rm(ii)} implies {\rm(i)}". The other half of Conjecture B,
namely ``{\rm(i)} implies {\rm(ii)}", remains open.

\section{Statement of the results}

For a $\rho'=\beta'+i\gamma'$ let of all ordinates of zeros of
$\zeta(s),$ $\gamma_c$ be the one for which $|\gamma_c-\gamma'|$ is
smallest (if there are more than one such zero of $\zeta(s),$ take
$\gamma_c$ to be the imaginary part of any one of them).

The following lemma is an immediate consequence of Lemmas 2 and 3
of~\cite{Zh}.

\begin{lemma}
\label{lem:ZhLemma} Assume RH. Let $\rho=\frac{1}{2}+i\gamma$ be a
simple zero of $\zeta(s)$ with $\gamma>0.$ Then
$$
\sum_{\beta'>\frac{1}{2}}\frac{\beta'-\frac{1}{2}}{(\beta'-\frac{1}{2})^2+(\gamma-\gamma')^2}=
\frac{1}{2}\log\gamma+O(1).
$$
\end{lemma}

Assuming RH and that $\frac{1}{2}+i\gamma_c$ is a simple zero, we
have
$$
\frac{\beta'-\frac{1}{2}}{(\beta'-\frac{1}{2})^2+(\gamma_c-\gamma')^2}\le
\frac{1}{2}\log\gamma_c+O(1).
$$
Hence, if $(\beta'-\tfrac{1}{2})\log\gamma'$ is small (which may
happen, since (i) is believed to be true), then
$|\gamma_c-\gamma'|\gg\sqrt{\frac{\beta'-\frac{1}{2}}{\log\gamma'}}.$
Our Theorem 1 may cause one to believe that $|\gamma_c-\gamma'|\ll
\sqrt{\frac{|\beta'-\tfrac{1}{2}|}{\log\gamma'}}$ for all
sufficiently large $\gamma'$. This may in turn suggest
$$
|\gamma_c-\gamma'|\asymp
\sqrt{\frac{|\beta'-\tfrac{1}{2}|}{\log\gamma'}},
$$
although one might also suspect that the right-hand side is off by a
factor of size a power of $\log\log\gamma'$, where the power may
vary depending on the size of $\beta' - \tfrac{1}{2}$ (the power may
become as high as $\tfrac{1}{2}$ for $\beta' - \tfrac{1}{2} \gg 1$)
in view of the conjecture made by Farmer, Gonek and
Hughes~\cite{FGH} based upon arguments from random matrix theory
that $\limsup\limits_{t\to\infty}\frac{|S(t)|}{\sqrt{\log t \log\log
t}}=\frac{1}{\pi\sqrt{2}}.$

\begin{theorem}
\label{thm:ordinate2} For any zero $\beta'+i\gamma'$ of $\zeta'(s)$
with a large $\gamma',$ there exists $\gamma_n$ such that
$\gamma'-1\le \gamma_n\le \gamma_{n+2}\le \gamma'+1$ and
$$
\min\{|\gamma_c-\gamma'|\log\gamma', \;
|\gamma_{n+2}-\gamma_n|\log\gamma_n\}\ll
\left(|\beta'-\tfrac{1}{2}|\log\gamma'\right)^{\tfrac{1}{2}}.
$$
\end{theorem}

Note that we haven't formulated the result in Theorem 1 in terms of
$\gamma_{n+1}-\gamma_{n}$ because if there are infinitely many zeta
zeros off the critical line, then since these zeros occur
symmetrically with respect to the critical line this difference will
be trivially $0$ infinitely often. In fact the statement of
Theorem~\ref{thm:ordinate2} holds more generally with
$\gamma_{n+n_0}$ in place of $\gamma_{n+2},$ where $n_0$ is any
fixed integer.

We also obtain unconditionally the following upper-bound.

\begin{theorem}
\label{thm:ordinate1} For any zero $\beta'+i\gamma'$ of $\zeta'(s)$
we have
$$
|\gamma_c-\gamma'|\ll |\beta'-\tfrac{1}{2}|^{\frac{1}{2}}.
$$
\end{theorem}

Besides the two statements in Conjecture B, let us pose the
following statement:
\begin{itemize}
\item[(iii)] \quad $\liminf\limits_{\gamma'\to\infty}|\gamma_c-\gamma'|\log\gamma'=0.$
\end{itemize}
In particular, from Theorem~\ref{thm:ordinate2} we immediately see
that if {\rm(i)} holds, then either {\rm(iii)} is true or
$\liminf\limits_{n\to\infty}(\gamma_{n+2}-\gamma_n)\log\gamma_n=0.$

Combining Theorem~\ref{thm:ordinate2} with Zhang's result which was
mentioned at the end of  \S 1 we derive

\begin{cor}
\label{cor:ordinate3} Assume that RH and {\rm(i)} hold. Then
{\rm(iii)} is true.
\end{cor}

Conjecture B claims that, under RH, (i) implies (ii). We establish
the following weaker result.

\begin{theorem}
\label{thm:ordinate4} Assume RH and
$\liminf\limits_{\gamma'\to\infty}
(\beta'-\frac{1}{2})(\log\gamma')(\log\log\gamma')^{2}=0$. Then
$$
\liminf\limits_{n\to\infty}(\gamma_{n+1}-\gamma_n)(\log\gamma_n)=0.
$$
\end{theorem}

We briefly recount some known conditional results related to
Theorems~\ref{thm:ordinate2} and~\ref{thm:ordinate1}. Guo~\cite{Guo}
(see also~\cite{Zh} for a generalization) has proved, under RH, if
for $\rho'$ with $T\le\gamma'\le 2T$ and
$\frac{1}{2}<\beta'<\frac{1}{2}+g(T)$ (where $g(T)\to 0$ as
$T\to\infty$) there exists a zero $\rho_1'=\beta_1'+i\gamma_1'$ of
$\zeta'(s)$ such that $|\rho_1'-\rho'|\ll\beta'-\frac{1}{2},$ then
$|\gamma_c-\gamma'|\ll \beta'-\tfrac{1}{2}.$ In the light of the
foregoing discussion, in Guo's result the condition of the existence
of such a zero $\rho_1'$ is crucial and probably can not be removed.
Zhang's paper contains the following result implicitly (see
(3.5)-(3.6) of~\cite{Zh}). Assume RH and
$\gamma_{n+1}-\gamma_n>\frac{2\pi\lambda}{\log T}$ with
$\gamma_n>\frac{T}{\log T},$ where $\lambda>1$ is such that the
condition $ \#\Bigl\{n: n<N(T), \gamma_{n+1}-\gamma_n
>\frac{2\pi\lambda}{\log T} \Bigr\}>c_0T\log T $ is satisfied with a
constant $c_0>0$ (from~\cite{CGGGH} this condition is known to hold
with $\lambda=1.33$). Then, there exists $\rho'$ such that
$|\rho'-\rho_n|<\frac{\nu}{\log T},$ where $\nu$ is such that
$(\frac{\nu}{\nu+2\pi\lambda})^2>\frac{\lambda+1}{2\lambda}.$

\section{Preliminaries}

We shall use some well-known properties of $\zeta(s)$ which can be
found in~\cite{KV} or ~\cite{T}. We recall the functional equation
\begin{equation}
\label{eqn:FE}
\pi^{-\frac{s}{2}}\Gamma(\tfrac{s}{2})\zeta(s)=\pi^{-\frac{1-s}{2}}\Gamma(\tfrac{1-s}{2})\zeta(1-s),
\end{equation}
and the partial fraction representation
$$
\frac{\zeta'(s)}{\zeta(s)}=b-\frac{1}{s-1}-\frac{1}{2}\frac{\Gamma'(\frac{s}{2}+1)}{\Gamma(\frac{s}{2}+1)}
+\sum_{\rho}\left(\frac{1}{s-\rho}+ \frac{1}{\rho}\right),
$$
where $b=-\frac{\gamma}{2}-1+\log
2\pi=-\sum\limits_{n=1}^{\infty}(\frac{1}{\rho_n}+\frac{1}{\overline{\rho_n}})+\frac{\log\pi}{2}.$
Using
\begin{equation}
\label{eqn:GammaLogDer} \frac{\Gamma'(s)}{\Gamma(s)}=\log |t|+O(1),
\qquad (0\le\sigma\le 2, |t|>2),
\end{equation}
we see that in the region $0\le\sigma\le 1, |t|>2$ the Riemann
zeta-function satisfies
$$
\frac{\zeta'(s)}{\zeta(s)}=\sum_{n=1}^{\infty}\left(\frac{1}{s-\rho_n}+\frac{1}{s-\overline{\rho_n}}\right)-
\tfrac{1}{2}\log |t| +O(1).
$$
Taking real parts and observing that in the region $0\le\sigma\le 1,
\, t>2$ the bound
$$
\sum_{n=1}^{\infty}{\Re\frac{1}{s-\overline\rho_n}}=\sum_{n=1}^{\infty}
{\frac{\sigma-\beta_n}{(\sigma-\beta_n)^2+(t+\gamma_n)^2}}=O(1)
$$
is valid, because  $\sum\limits_{n=1}^\infty\gamma_n^{-2}$ is
convergent and the $|\sigma-\beta_n|$ are bounded, we have
\begin{equation}
\label{eqn:ZetaLogDer} \Re
\frac{\zeta'(s)}{\zeta(s)}=\sum_{n=1}^{\infty}\frac{\sigma-\beta_n}{(\sigma-\beta_n)^2+(t-\gamma_n)^2}-
\tfrac{1}{2}\log t +O(1), \qquad (0\le\sigma\le 1, t>2).
\end{equation}
From the simple properties of the non-trivial zeros of $\zeta(s)$,
we know that
$$
\sum_{n=1}^{\infty}\frac{1}{1+(\gamma_n-T)^2}\ll \log T.
$$
for any real number $T\ge 2$. In particular, we know
\begin{equation}
\label{eqn:DistrZeros} \sum\limits_{|\gamma_n-T|\le 1}1\ll \log T,
\qquad \sum\limits_{|\gamma_n-T|\ge 1}\frac{1}{(\gamma_n-T)^2}\ll
\log T.
\end{equation}
It is also useful to remember that for every large $T>T_0>0,$
$\zeta(s)$ has a zero $\beta+i\gamma$ which satisfies
\begin{equation}
\label{eqn:SmallGapsOrdinates} |\gamma-T|\ll\frac{1}{\log\log\log
T},
\end{equation}
and that for any fixed $h,$ however small,
\begin{equation}
\label{eqn:SmallIntManyZeros} \sum_{{T\le \gamma_n\le T+h}}1>K\log
T, \qquad  (K=K(h)>0).
\end{equation}

The following lemma will play a role in the proof of
Theorems~\ref{thm:ordinate2} and~\ref{thm:ordinate1}.
\begin{lemma}
\label{lem:Lagrange} For any real numbers $a>0, \, x_1$ and $x_2,$
we have
$$
\left|\frac{x_1}{x_1^2+a}-\frac{x_2}{x_2^2+a}\right|\le
\frac{|x_1-x_2|}{a}.
$$
\end{lemma}
{\it Proof.} If $f(x)=\frac{x}{x^2+a},$ then
$$
|f'(x)|=\left|\frac{(x^2+a)-2x^2}{(x^2+a)^2}\right|\le\frac{x^2+a}{(x^2+a)^2}\le
\frac{1}{a},
$$
whence the result follows by the mean-value theorem.

\section{Proof of Theorems~\ref{thm:ordinate2} and~\ref{thm:ordinate1}}

We can assume that $|\beta'-\frac{1}{2}|$ is small, otherwise the
statements are trivial in view of~\eqref{eqn:SmallGapsOrdinates}
and~\eqref{eqn:SmallIntManyZeros}. We also can assume that
$\beta'\not= \frac{1}{2},$ since otherwise, $\beta'+i\gamma'$ is a
multiple zero of $\zeta(s)$ and again the results become trivial. We
also assume that $\gamma'$ is a large positive number and
$\gamma'\not=\gamma_n$ for any $n$.

Let $s=\sigma+it $ be in the region  $0\le \sigma\le 1, \, |t|>2.$
Taking logarithmic derivatives in the functional
equation~\eqref{eqn:FE}, and using~\eqref{eqn:GammaLogDer}, we have
$$
\frac{\zeta'(s)}{\zeta(s)}+\frac{\zeta'(1-s)}{\zeta(1-s)}=-\log
|t|+O(1).
$$
Since $\zeta'(\beta'+i\gamma')=\zeta'(\beta'-i\gamma')=0$, setting
$s=\beta'-i\gamma'$ we obtain
\begin{equation}
\label{eqn:differenceLogDer} \Re
\frac{\zeta'(\beta'+i\gamma')}{\zeta(\beta'+i\gamma')}-\Re
\frac{\zeta'(1-\beta'+i\gamma')}{\zeta(1-\beta'+i\gamma')}=\log\gamma'+O(1).
\end{equation}
Calculating the left-hand side of~\eqref{eqn:differenceLogDer}
via~\eqref{eqn:ZetaLogDer} with $s=\beta'+i\gamma'$ and
$s=1-\beta'+i\gamma'$ gives
$$
\sum_{n=1}^{\infty}\left({\frac{\beta'-\beta_n}{(\beta'-\beta_n)^2+(\gamma'-\gamma_n)^2}}-
{\frac{1-\beta'-\beta_n}{(1-\beta'-\beta_n)^2+(\gamma'-\gamma_n)^2}}\right)=\log\gamma'+O(1),
$$
so that
$$
\sum_{n=1}^{\infty}\left|{\frac{\beta'-\beta_n}{(\beta'-\beta_n)^2+(\gamma'-\gamma_n)^2}}-
{\frac{1-\beta'-\beta_n}{(1-\beta'-\beta_n)^2+(\gamma'-\gamma_n)^2}}\right|\ge\log\gamma'+O(1).
$$
Using Lemma~\ref{lem:Lagrange} with
$$
x_1=\beta'-\beta_n,\quad x_2=1-\beta'-\beta_n,\quad
a=(\gamma'-\gamma_n)^2>0,
$$
we have
\begin{eqnarray*}
&&\sum_{n=1}^{\infty}\left|{\frac{\beta'-\beta_n}{(\beta'-\beta_n)^2+(\gamma'-\gamma_n)^2}}-
{\frac{1-\beta'-\beta_n}{(1-\beta'-\beta_n)^2+(\gamma'-\gamma_n)^2}}\right|\\
&&\quad \le
2|\beta'-\tfrac{1}{2}|\sum_{n=1}^{\infty}\frac{1}{(\gamma'-\gamma_{n})^2}.
\end{eqnarray*}
Thus, we obtain
\begin{equation}
\label{eqn:Ready}
2|\beta'-\tfrac{1}{2}|\sum_{n=1}^{\infty}\frac{1}{(\gamma'-\gamma_n)^2}\ge\log\gamma'+O(1).
\end{equation}

First we prove Theorem~\ref{thm:ordinate2}. Without loss of
generality, we can assume that
$$
2|\beta'-\tfrac{1}{2}|\sum_{n=1}^{\infty}\frac{1}{(\gamma'-\gamma_{2n})^2}\ge
2|\beta'-\tfrac{1}{2}|{\sum_{n=1}^{\infty}}\frac{1}{(\gamma'-\gamma_{2n-1})^2}
$$
Then,
\begin{equation}
\label{eqn:Ready2}
4|\beta'-\tfrac{1}{2}|\sum_{n=1}^{\infty}\frac{1}{(\gamma'-\gamma_{2n})^2}\ge\log\gamma'+O(1).
\end{equation}
We recall that $\gamma'$ is a large number and
$|\beta'-\tfrac{1}{2}|$ is small. Then, from~\eqref{eqn:Ready2} we
have
$$
|\beta'-\tfrac{1}{2}|\sum_{n=1}^{\infty}\frac{1}{(\gamma'-\gamma_{2n})^2}\ge\frac{\log\gamma'}{8}.
$$
According to~\eqref{eqn:DistrZeros}, for a small
$|\beta'-\tfrac{1}{2}|,$ we have
$$
|\beta'-\tfrac{1}{2}|\sum_{|\gamma_n-\gamma'|\ge
1}\frac{1}{(\gamma'-\gamma_n)^2}\le \frac{\log\gamma'}{16},
$$
which implies
\begin{equation}
\label{readyTh2} |\beta'-\tfrac{1}{2}|\sum_{|\gamma_{2n}-\gamma'|\le
1}\frac{1}{(\gamma'-\gamma_{2n})^2}\ge\frac{\log\gamma'}{16}.
\end{equation}
Denote
$$
\delta=C\sqrt{\frac{|\beta'-\tfrac{1}{2}|}{\log\gamma'}},
$$
where $C>0$ will be chosen to be sufficiently large. Divide the
interval $[\gamma'-1, \gamma'+1]$ into small subintervals of the
type
$$
I_k=[\gamma'+k\delta, \gamma'+(k+1)\delta]\cap [\gamma'-1,
\gamma'+1],
$$
where $k$ runs through integers and $-\frac{1}{\delta}-1\le k\le
\frac{1}{\delta}+1.$ If there exists $n$ such that $ \gamma_{2n}\in
I_{-1}\cup I_0 $, then
$$
|\gamma_{2n}-\gamma'|\log \gamma'\le \delta\log
\gamma'=C\sqrt{|\beta'-\tfrac{1}{2}|\log\gamma'}
$$
and we are done in this case. Otherwise, we can
rewrite~\eqref{readyTh2} in the form
$$
|\beta'-\tfrac{1}{2}|\sum_{-\frac{1}{\delta}-1\le k\le -2}\quad
\sum_{\gamma_{2n}\in I_k}\frac{1}{(\gamma'-\gamma_{2n})^2}+
|\beta'-\tfrac{1}{2}|\sum_{1\le k\le \frac{1}{\delta}+1}\quad
\sum_{\gamma_{2n}\in
I_k}\frac{1}{(\gamma'-\gamma_{2n})^2}\ge\frac{\log\gamma'}{16}.
$$
If for some $k$ we have two numbers $n_1, n_2$ with
$\gamma_{2n_1}\in I_k,$ and $\gamma_{2n_2}\in I_k,$ then
$$
|\gamma_{2n_1}-\gamma_{2n_2}|\log\gamma'\le \delta\log
\gamma'=C\sqrt{|\beta'-\tfrac{1}{2}|\log\gamma'}
$$
and we are done in this case too. For this reason we can assume that
for a given $k$ we have at most one $n$ with $\gamma_{2n}\in I_k,,$
in which case we have
$$
|\gamma_{2n}-\gamma'|\ge |k+1|\delta, \quad {\rm when} \quad
-\frac{1}{\delta}-1\le k\le -2,
$$
and
$$
|\gamma_{2n}-\gamma'|\ge k\delta, \quad {\rm when} \quad 1\le k\le
\frac{1}{\delta}+1.
$$
Hence we have
$$
\frac{\log\gamma'}{16}\le
|\beta'-\tfrac{1}{2}|\sum_{|k|>0}\frac{1}{k^2\delta^2}\le
4C^{-2}\log\gamma'.
$$
However, this can not hold if $C>8.$ The proof of
Theorem~\ref{thm:ordinate2} is finished.

\bigskip

Now we proceed to prove Theorem~\ref{thm:ordinate1}. Let
$d=\min\limits_{n}|\gamma'-\gamma_n|.$ If we prove that $d\ll
|\beta'-\frac{1}{2}|^{\frac{1}{2}},$ then we are done.
From~\eqref{eqn:DistrZeros} we know
$$
2|\beta'-\tfrac{1}{2}|\sum_{d\le
|\gamma'-\gamma_n|<1}\frac{1}{(\gamma'-\gamma_n)^2}\ll\frac{|\beta'-\frac{1}{2}|\log\gamma'}{d^2},
$$
and
$$
2|\beta'-\tfrac{1}{2}|\sum_{|\gamma'-\gamma_n|\ge
1}\frac{1}{(\gamma'-\gamma_n)^2}\ll |\beta'-\tfrac{1}{2}|\log
\gamma'.
$$
Plugging these estimates into~\eqref{eqn:Ready} and recalling the
fact that $\gamma'$ is large, we obtain
$$
\frac{|\beta'-\frac{1}{2}|\log\gamma'}{d^2}
+|\beta'-\tfrac{1}{2}|\log\gamma'\gg \log\gamma'.
$$
Since $d=o(1)$ by~\eqref{eqn:SmallGapsOrdinates}, this reduces to
$$
\frac{|\beta'-\frac{1}{2}|\log\gamma'}{d^2}\gg\log\gamma',
$$
completing the proof of Theorem~\ref{thm:ordinate1}.

\section{Proof of Corollary~\ref{cor:ordinate3}}

In order to prove Corollary~\ref{cor:ordinate3}, assume that (iii)
is not true, i.e.
\begin{equation}
\label{eqn:corord2fromContradict}
\liminf_{\gamma'\to\infty}|\gamma_c-\gamma'|\log\gamma'\ge c_1>0.
\end{equation}
Then, there exists a constant $T_0$ such that
$$
|\gamma-\gamma'|\log\gamma' > \frac{c_1}{2}
$$
for $\gamma>T_0,\, \gamma'>T_0,$ and there can be at most finitely
many multiple zeros of $\zeta(s)$. Hence, assuming RH,
Theorem~\ref{thm:ordinate2} implies
$$
\liminf_{\gamma\to\infty}(\gamma^+-\gamma)\log\gamma=0.
$$
Therefore, in Zhang's result which was mentioned at the end of \S 1,
we can take $\alpha_1$ to be small (and therefore, we can take
$\alpha_2$ to be small too) and deduce that for any $\rho$ with a
large $\gamma$ there exists a zero $\rho'$ of $\zeta'(s)$ such that
$$
|\gamma-\gamma'|\log\gamma\le |\rho'-\rho|\log\gamma<\frac{c_1}{2}.
$$
This contradicts~\eqref{eqn:corord2fromContradict} and proves
Corollary~\ref{cor:ordinate3}.

\section{Proof of Theorem~\ref{thm:ordinate4}}

The proof presented here stems from an idea of Haseo Ki. We now work
under the assumptions that the RH is true, and
\begin{equation}
\label{eqn:hypologlog2} \liminf_{\gamma' \to\infty}(\beta' -
\tfrac{1}{2})(\log\gamma') (\log\log\gamma')^{2} = 0 .
\end{equation}
For our purpose we may also assume that all but finitely many of the
zeta zeros are simple, because otherwise
$\lim\inf(\gamma_{n+1}-\gamma_n)\log\gamma_n=0 $ holds trivially.

For a $\rho' = \beta' + i\gamma'$, member of a sequence with the
property ~\eqref{eqn:hypologlog2}, there are two possibilities:
\newline Either
$$
|\gamma_{c} - \gamma'| \leq |\gamma_{n+1} - \gamma_{n}| , \quad
\forall \gamma_{n}\in [\gamma' - 1, \gamma' +1],
$$
or
$$
\exists \gamma_{n}\in [\gamma' - 1, \gamma' +1], \quad |\gamma_{n+1}
- \gamma_{n}| < |\gamma_{c} - \gamma'|.
$$

If there is a subsequence of $\rho'$ satisfying the second
possibility, we have by Theorem 1 for the corresponding $\gamma_n$,
$$
|\gamma_{n+1} - \gamma_{n}|(\log \gamma_{n}) \ll (|\beta' -
\tfrac{1}{2}|\log \gamma')^{\tfrac{1}{2}}.
$$
Thus, in this case we don't even need the full strength of the
condition ~\eqref{eqn:hypologlog2} to conclude that
$\displaystyle\liminf_{n\to\infty}(\gamma_{n+1}-\gamma_{n})\log\gamma_n
= 0$.

From now on we may take that after a point on all $\rho'$ from a
sequence with the property ~\eqref{eqn:hypologlog2} satisfy the
first possibility. Suppose
$$
\lim\inf_{n\to\infty}(\gamma_{n+1}-\gamma_n)\log\gamma_n>0 ,
$$
so that there exists a fixed $\delta>0$ such that
$$
\gamma_{n+1}-\gamma_n>\frac{\delta}{\log\gamma_n}
$$
for all sufficiently large $n$.

We apply the formula ~\cite[Theorem 9.6 (A)]{T}
$$
\frac{\zeta'}{\zeta}(s)= \sum_{|\gamma-t|\leq
1}\frac{1}{s-\rho}+O(\log t)
$$
at $s=\rho'$, where $\rho'$ is a member of a sequence obeying
~\eqref{eqn:hypologlog2}. So we can write
\begin{equation}
\label{eqn:96Arho'} 0=\frac{1}{\rho'-\rho_c}+\sum_{\rho\neq\rho_c
\atop |\gamma-\gamma'|\leq 1 }\frac{1}{\rho'-\rho}+O(\log\gamma').
\end{equation}
We now examine the sum occuring in this formula. Clearly,
$$
| \sum_{\rho\neq\rho_c \atop |\gamma-\gamma'|\leq 1}
\frac{1}{\rho'-\rho} | \leq \sum_{\gamma\neq\gamma_c \atop
|\gamma-\gamma'|\leq 1} \frac{1}{|\gamma-\gamma'|}.
$$
By our assumption we have, for all positive integers $j$,
$$
|\gamma_{c\pm j} -\gamma'| \geq {j\delta\over 3\log\gamma'}
$$
(here $\gamma_{c\pm j}= \gamma_{n_{0}\pm j}$ when
$\gamma_{c}=\gamma_{n_{0}}$). Since the sum is over the zeros with
$\gamma$ in an interval of radius $1$ around $\gamma'$, we see that
$j$ can be at most as large as ${\kappa\log\gamma'\over \delta}$
with some absolute constant $\kappa$. Therefore
$$
\sum_{\gamma\neq\gamma_c \atop |\gamma-\gamma'|\leq 1}
\frac{1}{|\gamma-\gamma'|} \ll {(\log\gamma')(\log\log\gamma')\over
\delta}.
$$
Hence we can rewrite ~\eqref{eqn:96Arho'} as
$$
0=\frac{1}{\rho'-\rho_c}+
O\Bigl({(\log\gamma')(\log\log\gamma')\over \delta}\Bigr),
$$
from which we see that
\begin{equation}
\label{eqn:bbound}
\frac{1}{\sqrt{(\beta'-\tfrac{1}{2})^2+(\gamma'-\gamma_c)^2}} \leq
{\kappa_{1}(\log\gamma')(\log\log\gamma')\over \delta}
\end{equation}
for some absolute constant $\kappa_{1}$. Now recall that $\rho'$
satisfies the first possibility, so that by Theorem 1 we have
$$
|\gamma' - \gamma_c| \leq
\kappa_{2}\Bigl(\frac{\beta'-\tfrac{1}{2}}{\log\gamma'}\Bigr)^{\tfrac{1}{2}}
$$
for some absolute constant $\kappa_{2}$. Using this in
~\eqref{eqn:bbound} we get
$$
\frac{1}{(\beta'-\tfrac{1}{2})^2+ \kappa_{2}^{2}
(\frac{\beta'-\tfrac{1}{2}}{\log\gamma'}) } \leq
{\kappa_{1}^{2}(\log\gamma')^2 (\log\log\gamma')^2 \over \delta^2}.
$$
Now the quadratic formula yields
$$ \beta'-\tfrac{1}{2} \geq \frac{\delta^2}{2(\kappa_{1}\kappa_{2})^{2}
(\log\gamma')(\log\log\gamma')^2}
$$
for sufficiently large $\gamma'$, which contradicts the assumption
~\eqref{eqn:hypologlog2}.

This completes the proof of Theorem~\ref{thm:ordinate4}.

\bigskip

{\bf Acknowledgement.} This work was supported by Project
PAPIIT-IN105605 from the UNAM. Y{\i}ld{\i}r{\i}m is grateful to
Instituto de Matem\'aticas, UNAM, Campus Morelia, M\'exico for its
hospitality. The authors are grateful to Professor Haseo Ki for
sharing his thoughts which led to the present form of Theorem 3,
improving upon the former version greatly.

\end{document}